\documentclass[11pt]{article}

\usepackage{cite}                                       
\usepackage{graphicx}        
\usepackage[latin1]{inputenc}
\usepackage[english]{babel}
\usepackage{amssymb,amsmath,amsfonts}
\usepackage{amsthm,amscd}
\usepackage{fancyhdr}

\def\today{25.05.09}

\def\eps{\epsilon}
\def\im{\mathrm{i}}
\def\sech{{\rm sech}}
\def\E{\mathcal{E}}
\def\C{\mathcal{C}}

\def\S{\mathcal{S}}
\def\L{\mathcal{L}}

\def\U{\mathcal{U}}

\def\RR{\mathbb{R}}
\def\ZZ{\mathbb{Z}}
\def\Q{\mathcal{Q}}

\def\tond#1{\left(#1\right)}
\def\quadr#1{\left[#1\right]}
\def\graff#1{\left\{#1\right\}}
\def\norma#1{\left\| #1\right\|}
\def\inter#1{\langle#1\rangle}

\newcommand{\Sc}{{\mathcal S}}
\newcommand{\Hc}{{\mathcal H}}

\newcommand{\Ic}{{\mathcal I}}

\newcommand\nqmu[1]{\left\|#1\right\|_{\Q_\mu}}
\newcommand\nll[1]{\left\|#1\right\|_{\ell^2}}

\newtheorem{theorem}{Theorem}[section]
\newtheorem{lemma}{Lemma}[section]

\newtheorem{proposition}{Proposition}[section]
\newtheorem{definition}{Definition}[section]
\newtheorem{remark}{Remark}[section]

\title{Continuous approximation of breathers \\ in one and two
  dimensional DNLS lattices.}

\author{D.~Bambusi and T.~Penati}
\date{\today}

\begin{document}

\maketitle

\begin{abstract}
In this paper we construct and approximate breathers in the DNLS model
starting from the continuous limit: such periodic solutions are
obtained as perturbations of the ground state of the NLS model in
$H^1(\RR^n)$, with $n=1,2$. In both the dimensions we recover the
Sievers-Takeno (ST) and the Page (P) modes; furthermore, in $\RR^2$
also the two hybrid (H) modes are constructed. The proof is based on
the interpolation of the lattice using the Finite Element Method
(FEM).
\end{abstract}

\section{Introduction}
\label{intro}

In this paper we study the problem of constructing breathers in the
one and two dimensional discrete nonlinear Schr\"odinger (DNLS)
equation starting from the continuous limit.

The breathers we construct are critical points of the Hamiltonian
function constrained to the surface of constant $\ell^2$ norm. Such
critical points are obtained by continuation from the continuous model
constituted by the nonlinear Schr\"odinger (NLS) equation. The
connection between the discrete and the continuous system is obtained
by using the finite elements (FEM). {\it This allows to identify the
  phase space of the discrete system with a subspace of the phase
  space of the continuous system. }

\begin{figure}[t]
\begin{center}
\includegraphics[width=.8 \textwidth]{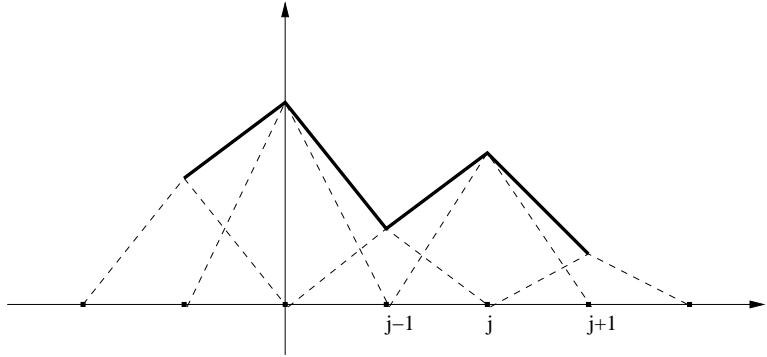}
\end{center}
\caption{Finite element basis $s_j(x)$ and interpolating function
  $\psi(x)$.}
\label{fig1}
\end{figure}

For example, consider the one dimensional case. The space of the
finite elements is constructed as follows: first we associates to the
$j$-th point of the discrete lattice a continuous piecewise linear
function $s_j(x)$, whose value is $1$ at $x=j$ and which vanishes for
$|x - j|\geq 1$ (see Fig.~\ref{fig1}).  To a sequence $\psi_j$, we
associate the function $\psi(x):=\sum_j\psi_js_j(x/\mu)$, where
$\mu>0$ is a small parameter representing the mesh of the lattice. The
space generated by the functions $s_j(x/\mu)$ will be denoted by
$\E_\mu$.

Once this is done one can compare the functionals of the continuous
system and those of the discrete one. In order to do this, denote by
$H_c$ and $N_c$ the Hamiltonian and the norm of the continuous system,
and consider the restriction of such functions to the space of the
discrete system $\E_\mu$. By the standard theory of integration one
can say that the restricted functionals are close to the Hamiltonian
$H_d$ and the norm $N_d$ of the discrete system. So the idea is to
consider a non-degenerate critical point of the functional of the
continuous system, a critical point laying close to the manifold of
the finite elements and to continue such a critical point to a
critical points of the discrete functional.

However there is a delicate point in the game: namely that the
difference between the discrete functional and the continuous one
should be small {\it when the phase space is endowed with the energy
  norm}. This turns out to be true thanks to a special property of the
finite elements: the fact that one has
\begin{displaymath}
\int_{\RR^n}|\nabla \psi|^2dx = \mu^{n}\sum_{|{\bf
    j-l}|=1}\frac{|\psi_{\bf j} - \psi_{\bf l}|^2}{\mu^{2}},
\end{displaymath}
with no error. Due to this property the difference between the
continuous and the discrete functional turns out to be a functional
which is {\em small and smooth} on the energy space. This allows to
apply the {\em implicit function theorem} and to continue critical
points of the continuous system to critical points of the discrete
one.

In order to be concrete we study in detail a one dimensional and a two
dimensional model. We use known results on existence and
non-degeneracy of the ground state of the continuous system in order
to apply the above theory. In these paper we construct two
(resp. four) kinds of discrete breathers in the 1-(resp. 2)
dimensional case, which are the continuation of the continuous
breather. In order to avoid problems related to the translational
invariance of the continuous system we work here in spaces of
reflection invariant sequences. Thus the breathers we find for the
discrete system are reflection invariant too.

In dimension one, the breather of the first kind is centered at a
lattice site and corresponds to the so called Sievers-Takeno mode
(ST), while the breather of the second kind is centered in the middle
of a cell of the lattice and corresponds to the so called Page mode
(P). In dimension two, besides the ST and P modes, we have two other
localized solutions, usually called hybrid (H) modes since they are
centered in the middle of one of the two face of the cell.

As far as we know, the result of the present paper is the first one in
which the continuous approximation is used in order to construct exact
breathers of a lattice model. In dimension 1 the method of spatial
dynamics also allows to construct and approximate breathers (see
\cite{Jam03}). However such a method is strictly one dimensional,
while our method in principle applies to any dimension. Existence of
breathers was also proved variationally in \cite{Wei_MI99} and in
\cite{AubKopKad01}, but such methods do not allow to approximate the
breathers and only allow to find one breather for each
model. Breathers in DNLS have also been widely studied numerically
(see for example \cite{KevrRasBis01B,JmsKevr08,FlaWil98}.

The main advantage of our method is that it is quite flexible and
allows to directly deduce informations on the shape of the breather
starting from the continuous limit. 

We recall that the possibility of using the continuous limit in order
to approximate the dynamics of discrete systems has been widely
investigated, in particular we recall the papers
\cite{BCP,Sch98,KSM92,SW99,BP06,BCP09} in which an approximation valid
for long but finite times and the papers
\cite{FPI,FPII,FPIII,FPIV,HofWay08,MizPeg08} where an infinite time
approximation has been obtained.

The plan of the paper is the following. In Section \ref{Main} we
present the result and motivate our continuum limit approach. In
Section \ref{IFT} we formulate in Theorem \ref{t.a.1} the Implicit
Function Theorem applied to our problem and in Section \ref{breath} we
construct the FEM to interpolate the discrete model and we verify the
hypothesis of Theorem \ref{t.a.1}.


\section{Main result.}
\label{Main}

We study here the discrete focusing nonlinear Schr\"odinger equation
(DNLS) in $\RR^n$ with $n=1,2$
\begin{equation}
\label{DNLS}
\im \dot \psi_{\bf l} = -
\frac1{\mu^{2}}(\Delta_1\psi)_{\bf l} - |\psi_{\bf l}|^{2p}\psi_{\bf
  l}\ ,\qquad {\bf l}\in\ZZ^n,\quad \frac12\leq p<\frac2n
\end{equation}
where $\Delta_1$ is the $n$-dimensional discrete Laplacian defined by
\begin{eqnarray*}
(\Delta_1\psi)_j &:=& \psi_{j+1} + \psi_{j-1}-2\psi_j,
  \\ (\Delta_1\psi)_{j,k} &:=& (\psi_{j+1,k} +
  \psi_{j-1,k}-2\psi_{j,k}) + (\psi_{j,k+1} +
  \psi_{j,k-1}-2\psi_{j,k}),
\end{eqnarray*}
and $\mu$ is the {\em lattice mesh}. In particular we look for
solutions of the form
\begin{equation}
\label{e.1}
\psi_{\bf l}(t)=e^{-\im \lambda t}\psi_{\bf l}.
\end{equation}
Then the sequence $\psi_{\bf l}$ fulfils
\begin{equation}
\label{lagrange}
\lambda \psi_{\bf l} = -\frac1{\mu^2} (\Delta_1\psi)_{\bf l} -
|\psi_{\bf l}|^{2p}\psi_{\bf l},
\end{equation}
and thus it is a {\em critical} point of the Hamiltonian function
\begin{equation}
\label{e.2}
H_d:= \mu^n\quadr{\frac12 \sum_{|{\bf j-l}|=1}\frac{|\psi_{\bf j} - \psi_{\bf
    l}|^2}{\mu^{2}} - \frac1{p+1}\sum_{{\bf l}\in\ZZ^n}|\psi_{\bf
  l}|^{2p+2}}
\end{equation}
constrained to a surface of constant value of the norm
\begin{equation}
\label{e.3}
N_d:=\mu^n\sum_{{\bf l}\in\ZZ^n}|\psi_{\bf l}|^2,
\end{equation}
where the factors $\mu^n$ have been inserted for future convenience.
The main result of the present paper consists in showing that such a
solution can be {\bf constructed and approximated} starting from the
continuous model constituted by the Nonlinear Schr\"odinger Equation
(NLS), namely
\begin{equation}
\label{NLS}
\im \dot \psi=-\Delta\psi-|\psi|^{2p}\psi.
\end{equation}
More precisely, consider the Hamiltonian $H_c$ and (the square of) the
$L^2$ norm $N_c$, given by
\begin{equation}
\label{e.31}
H_c:=\int_{\RR^n}\left[|\nabla\psi|^2-\frac{1}{p+1}|\psi|^{2p+2}\right]
,\qquad N_c:=\int_{\RR^n}|\psi|^2,
\end{equation}
then a periodic solution $\psi(x,t)=e^{-\im \lambda t}\psi(x)$ of
\eqref{NLS} fulfils the following continuous approximation of
\eqref{lagrange}
\begin{equation}
\label{lagrange.c}
\lambda\psi = -\Delta\psi-|\psi|^{2p}\psi.
\end{equation}
According to classical results on \eqref{lagrange.c} (see
\cite{BeLio83, BeLioPel81,ColGlaMar78}), there exists a unique real
valued, positive, radially symmetric and exponentially decaying
function $\psi_c$ which realizes the minimum of
$H_c\vert_{N_c=1}$. For example, in the case $n=1$ and $p=1$ it can be
computed explicitly
\begin{equation}
\label{1D.sol}
\psi_c(x):=\frac{1}{\sqrt 2}\sech\left( \frac{x}{2}\right).
\end{equation}
If we interpret the discrete functionals $H_d,\,N_d$ as
$\mu$-perturbations of $H_c,\,N_c$ and we restrict to a class of
``even'' functions in order to remove any possible degeneracy of the
minimum $\psi_c$, then we can continue the solution $\psi_c$ of
\eqref{lagrange.c} to a solution $\psi(\mu)$ of \eqref{lagrange}. 

In order to state the precise result we are going to prove, we first
need to define the configuration space $\Q_\mu$ for $\psi_{\bf l}$:

\begin{definition}
\label{d.Qspace}
The space $\ell^2(\ZZ^n,\RR)$ will be denoted by $\Q_\mu$ when endowed
with the norm
\begin{equation}
\label{e.Qnorm}
  \nqmu{\psi}^2 := \mu^n\nll{\psi}^2 +
  \frac1{\mu^{2-n}}\inter{\psi,-\Delta_1 \psi}_{\ell^2}^2.
\end{equation}
\end{definition}

\begin{theorem}
\label{main.theo}
For any $\mu$ small enough and $\frac12\leq p<\frac2n$ there exist
$2^n$ distinct real valued sequences $\psi_{\bf l}^i(\mu)$ which are
solutions of \eqref{lagrange}. Such solutions are even sequences
$\psi_{-\bf l} = \psi_{\bf l}$ lying on the surface $N_d=1$. One has
\begin{equation}
\label{est}
\nqmu{\psi^i - \Psi^i}\leq C_1\mu,\qquad \sup_{{\bf l}\in\ZZ^n}
\left|{\psi^i_{\bf l}(t) - \Psi^i_{\bf l}(t)} \right| \leq
C_2\mu^{\frac32-\frac{n}2}
\end{equation}
where $\Psi^i$ is defined by
\begin{eqnarray*}
&&\begin{cases} \Psi^1_{\bf l}:=
\psi_c\tond{\mu j},\\ \Psi^2_{\bf l}:=\psi_c\tond{\mu j +
  \frac\mu2},
\end{cases}\,n=1,\\
&&\begin{cases} \Psi^1_{\bf l}:=\psi_c\tond{\mu j, \mu
      k},\\ \Psi^2_{\bf l}:=\psi_c\tond{\mu j, \mu k +
      \frac\mu2},\\ \Psi^3_{\bf l}:=\psi_c\tond{\mu j + \frac\mu2,
      \mu k},\\ \Psi^4_{\bf l}:=\psi_c\tond{\mu j + \frac\mu2, \mu
      k + \frac\mu2},
\end{cases} \, n=2.
\end{eqnarray*}
\end{theorem}

\subsection{Comments.}

\begin{enumerate}
\item the first of \eqref{est} is not empty since by using
  \eqref{e.approx} we get
\begin{displaymath}
\nqmu{\Psi^i}\sim 1.
\end{displaymath}
Moreover, we stress that by its definition the approximating sequence
$\Psi^i_{\bf l}$ is bounded uniformly in $\mu$
\begin{displaymath}
\Psi^i_{\bf l}\leq \norma{\psi_c}_{L^\infty}
\end{displaymath}
but is localized\footnote{We can fix the set where $\Psi^i_{\bf l}$ is
  localized as $\Omega:=\{{\bf l}\in\ZZ^n \,s.t.\,\Psi^i_{\bf l}\geq
  \frac12\Psi_0\}$} on an increasing interval $[-k,k]^n$ with
$k\sim1/\mu$.

\item the first of \eqref{est} immediatly implies
\begin{displaymath}
\nll{\psi^i - \Psi^i}\leq C\mu^{1-\frac{n}2} \quad\Longrightarrow\quad
\left|{\psi^i_{\bf l}(t) - \Psi^i_{\bf l}(t)} \right| \leq
C_2\mu^{1-\frac{n}2}\qquad \forall {\bf l},
\end{displaymath}
an estimate which is empty in the case $n=2$. Lemma \ref{l.sup} in
Section 4.3 is necessary to improve the above result. We do not know
whether the exponent $\frac32-\frac{n}2$ is optimal or not.


\item We stress that the problem \eqref{lagrange} is equivalent to the
$\mu$-indepen\-dent one
\begin{equation}
\label{lagrangenomu}
\tilde\lambda \varphi_{\bf l} = - (\Delta_1\varphi)_{\bf l} -
|\varphi_{\bf l}|^{2p}\varphi_{\bf l}
\end{equation}
with the constrain
\begin{displaymath}
\sum_{{\bf l}\in\ZZ^n}{|\varphi_{\bf l}|^2} = E,\qquad E\ll 1.
\end{displaymath}
This can be seen by the scaling
\begin{displaymath}
\varphi_{\bf l} = \mu^{\frac1p} \psi_{\bf l},\qquad \tilde\lambda =
\mu^2\lambda,\qquad \mu^{\frac2p-n} = E,
\end{displaymath}
and observing that
\begin{displaymath}
\frac2p-n>0 \qquad\Longleftrightarrow\qquad p<\frac2n.
\end{displaymath}
\end{enumerate}


\bigskip
\bigskip


\section{The Implicit Function Theorem.}
\label{IFT}

The situation we will meet is summarized in the following abstract
scheme. Let $\Hc$ be a Hilbert space, and for any $\mu$, let $\E_\mu$
be a subspace of $\Hc$. Let $H_c\in C^2(\Hc)$ and $N_c\in
C^\infty(\Hc)$ be two functionals, with $N_c$ being a
submersion. Correspondingly we define
\begin{displaymath}
\Sc:=\left\{ \psi\in \Hc\ :\ N_c(\psi)=1\right\}\ .
\end{displaymath}
Then we define the ``discrete'' objects: let
$H_d:=H_{\epsilon_1,\mu}\in C^2(\E_\mu)$ and
$N_d:=N_{\epsilon_2,\mu}\in C^\infty (\E_\mu)$ be functionals
depending smoothly on two additional parameters $\epsilon_1, \epsilon
_2$. Define
$$ \Sc_{\epsilon_2,\mu}:=\left\{ \psi\in
\E_\mu\ :\ N_{\epsilon_2,\mu}(\psi)=1\right\}\ .
$$
We make some assumptions.
\begin{itemize}
\item[i.]  There exists $\psi_c\in\Hc$ which is a coercive minimum of
  $H_c\big|_{\Sc}$, namely it is a minimum and fulfills
\begin{equation}
\label{coerc}
d^2H_c\big|_{\Sc}(\psi_c)(h,h)\geq C\norma{h}^2,\qquad h\in
T_{\psi_c}\Sc;
\end{equation}
moreover
\begin{equation}
\label{e.a.2}
d(\E_\mu,\psi_c)\leq C\mu\ ,
\end{equation}
for all $\mu$ small enough.
\end{itemize}

Let $\psi_0\in \E_\mu$ be such that $\norma{\psi_c-\psi_0}\leq C \mu$
and let $\U\subset \E_\mu$ be an open neighborhood of $\psi_0$ then we
assume
\begin{itemize}
\item[ii.]
\begin{eqnarray}
\label{hp.ii}
  \norma{H_{\epsilon_1,\mu} - H_{0,\mu}}_{C^2(\U)}&\leq& C
  (\epsilon_1+\mu),\qquad H_{0,\mu}:= H_c\big|_{\E_\mu}
  \\ \norma{N_{\epsilon_2,\mu}-N_{0,\mu}}_{C^k(\U)}&\leq& C
  (\epsilon_2+\mu),\qquad N_{0,\mu}:= N_c\big|_{\E_\mu} \nonumber
\end{eqnarray}
for some large enough $k$
\end{itemize}

\begin{theorem}
\label{t.a.1}
Under the above assumptions, for any $\epsilon_1,\epsilon_2,\mu$ small
enough, there exists a unique $\psi_{\epsilon_1,\epsilon_2,\mu}$,
which is a coercive minimum of $H_d\big|_{N_d=1}$. Moreover one has
\begin{equation}
\label{dist}
\norma{\psi_{\epsilon_1,\epsilon_2,\mu} - \psi_c}\leq C'(\mu +
\epsilon_1+\epsilon_2)\ .
\end{equation}
\end{theorem}

\proof The result is local, so we restrict to a neighborhood of
$\psi_c$. Define 
\begin{eqnarray}
\Sc_{0,\mu} := \left\{ \psi\in \E_\mu\ :\ N_{0,\mu}(\psi)=1\right\} =
\Sc\cap\E_\mu
\end{eqnarray}
and take $\psi_0\in\S_{0,\mu}$. Remark that, due to smoothness of
$H_c$ one has
\begin{equation}
\label{e.dim.1}
\norma{d(H_c\big|_{\Sc_{0,\mu}})(\psi_0)}\leq C \mu\ .
\end{equation} 
By coercivity \eqref{coerc} and Lax-Milgam Lemma, the second
differential 
\begin{displaymath}
d^2(H_c\big|_{\Sc_{0,\mu}})(\psi_0):T_{\psi_0}\Sc_{0,\mu} \to
T^*_{\psi_0}\Sc_{0,\mu}
\end{displaymath}
defines an isomorphism bounded together with its inverse uniformly
with respect to all the parameters.

From assumption ii, there exists a local isomorphism
\begin{displaymath}
\Ic_{\epsilon_2,\mu}:\Sc_{0,\mu}\to \Sc_{\epsilon_2,\mu},
\end{displaymath}
which satisfies
\begin{displaymath}
\norma{\Ic_{\epsilon_2,\mu}-Id }_{C^k}\leq C(\epsilon_2+\mu)\ .
\end{displaymath}
The statement is then equivalent to the existence of a coercive
minimum of $H_{\eps_1,\mu}\circ \Ic_{\epsilon_2,\mu}$. To get it
remark that
\begin{equation}
\label{e.dim.2}
H_{\eps_1,\mu}\circ
\Ic_{\epsilon_2,\mu}=H_c\big|_{\Sc_{0,\mu}}+O(\epsilon_1+\epsilon_2+\mu)\ .
\end{equation}
Due to \eqref{e.dim.1} and \eqref{coerc}, the Implicit Function
Theorem applies and gives the result. \qed

%

\bigskip


\section{Applications to breathers}\label{breath}
In order to avoid gauge and the translational invariance of the
problem, in particular of the continuous system, we will work in a
space of real valued functions ``invariant'' under the involution
\begin{equation}
\label{S}
S:\psi(z)\mapsto\psi(-z),\qquad z\in\RR^n .
\end{equation}
More precisely, in $H^1(\RR^2,\RR)$ we will consider functions
fulfilling
\begin{equation}
\label{constr}
\int_{\RR^n}|\psi(z)-\psi(-z)|^2 = \norma{\psi - S\psi}^2_{L^2} = 0,
\end{equation}
which is equivalent to \eqref{S} almost everywhere and is a condition
well defined in $H^1(\RR^2,\RR)$. 
\begin{lemma}
\label{NDegen}
Let $\psi_c$ be a solution of \eqref{lagrange.c} with $p<\frac2n$ and
let
\begin{displaymath}
\Hc = \graff{\psi\in H^1(\RR^n,\RR): \norma{\psi - S\psi}^2_{L^2} =
  0} =: H^1_s, 
\end{displaymath}
then assumption \eqref{coerc} of Theorem \ref{t.a.1} holds.
\end{lemma}

\proof This Lemma directly follows from Proposition D.1 of
\cite{FGJS04} by remarking that $T_{\psi_c}\Sc\subset X$, with $X$
defined in the statement of Prop. D.1.\qed.

\bigskip

\begin{remark}
We stress that the constrain \eqref{constr} is {``natural''} for the
problem \eqref{lagrange}, since \eqref{S} is a symmetry for both the
Hamiltonian \eqref{e.2} and the Norm \eqref{e.3}. Hence, a critical
point for the restricted problem is also a critical point for the
original problem.
\end{remark}

\bigskip

In the following subsections we construct the linear manifold $\E_\mu$
of the finite elements, and prove the estimates \eqref{e.a.2} and
\eqref{hp.ii} for the two considered applications. We deal with the
ST-breather, since the other ones follow by small changes in the
definition of $\E_\mu$.


\subsection{The case $n=1$}

Let ${\bf l} = j$ and define the sequence of functions $s_j(x)$ by
\begin{equation}
\label{fe.1}
s_j(x)=
\begin{cases}
0,\qquad\qquad  {\rm if}\quad |x-j|>1\\
x-j+1,\quad  {\rm if}\quad   j-1\leq x\leq j\\
- x+j+1,\quad  {\rm if}\quad  j\leq x\leq j+1
\end{cases}
\end{equation}
and, to a sequence $\psi_j\in\Q_\mu$, we associate a function
\begin{equation}
\label{e.5}
\Psi(x) := \sum_j\psi_js_j(x/\mu).
\end{equation}
On the interval $T_j:=[\mu j,\mu(j+1))$ the above function reads
\begin{equation}
\label{e.interpol}
\Psi(x) = (x-\mu j)\frac{(\psi_{j+1} - \psi_{j})}{\mu} + \psi_{j}.
\end{equation}

\begin{definition}
\label{linspace.1}
We denote by $\E_\mu$ the linear space composed by the functions of
the form \eqref{e.5} with $\psi_j\in\Q_\mu$.
\end{definition}

The following Lemma gives the equivalence between {\bf the function
  space} $\E_\mu$ and {\bf the sequence space} $\Q_\mu$.

\begin{lemma}
\label{l.grad.1}
Let $\Psi\in\E_\mu$ then
\begin{equation}
\label{e.norm.der.1}
\int_{\RR}\Psi_x^2 = \frac1\mu\sum_{j\in\ZZ}{(\psi_{j+1} - \psi_{j})^2}.
\end{equation}
Moreover
\begin{equation}
\label{e.norm.ell2.1}
\mu\sum_{j}{\psi_{j}^2} = \int_{\RR} \Psi^2 +
\frac\mu3\int_{\RR}\Psi_x^2dx.
\end{equation}
\end{lemma}

\proof  Let us first decompose $\RR = \cup_{j\in\ZZ} T_j$. The weak
derivative of $\Psi$ is
\begin{displaymath}
\Psi_x(x) = \sum_{j\in\ZZ}{\frac{\tond{\psi_{j+1}-\psi_j}}\mu\chi_{(\mu
    j,\mu (j +1))}(x)}
\end{displaymath}
which gives immediately
\begin{equation}
\Psi^2_x(x) =
\quadr{\sum_{j\in\ZZ}\frac{(\psi_{j+1}-\psi_j)}\mu\chi_{T_j}(x)}^2 =
\sum_{j\in\ZZ}{\frac{(\psi_{j+1}-\psi_j)^2}{\mu^2}\chi_{T_j}(x)},
\end{equation}
since $\chi_{T_j}(x)\chi_{T_i}(x)=0$ for $i\not=j$. From
$\int_\RR{\chi_{T_j}(x)dx} =\mu$ one gets
\begin{equation}
\label{norm.der.1}
\int_\RR{\Psi_x^2dx} = \mu
\sum_{j\in\ZZ}{\frac{(\psi_{j+1}-\psi_j)^2}{\mu^2}} =
\frac1\mu\sum_{j\in\ZZ}{(\psi_{j+1} - \psi_{j})^2}.
\end{equation}
If we plug \eqref{e.interpol} in the integral $\norma{\Psi}^2_{L^2}$,
a direct computation gives the estimate \eqref{e.norm.ell2.1}.\qed

\bigskip
\bigskip

%
%

\begin{proposition}
\label{p.1d}
Let $\Psi\in\E_\mu$ be as in \eqref{e.5} and let us define
\begin{eqnarray*}
G_c(\Psi) &:=& \int_{\RR}|\Psi|^{q+2},\\
G_d(\Psi) &:=& \mu\sum_{j\in\ZZ}{|\psi_{j}|^{q+2}},\\
R_G(\Psi) &:=& G_c(\Psi) - G_d(\Psi);
\end{eqnarray*}
if $q\geq 1$ then $R_G\in\C^2(\E_\mu)$ and for any bounded open set
$\U\subset\E_\mu$, there exists $C(\U)$ such that
\begin{equation}
\label{e.p.A}
\norma{R_G}_{C^2(\U)}\leq C\mu.
\end{equation}
\end{proposition}

\proof The term $R_G$ can be represented through the Euler-\-MacLaurin
formula
\begin{equation}
\sum_{j\in\ZZ}{f(j)} = \int_{\RR}f(y)dy +
\int_{\RR}f_y(y)P_1(y)dy,\qquad P_1(s) = s-[s]-\frac12.
\end{equation}
Indeed, if we set $f(y) = |\Psi(\mu y)|^{q+2}$, we have
\begin{eqnarray*}
\mu\sum_{j\in\ZZ}{|\Psi_j|^{q+2}} &=& \mu\sum_{j\in\ZZ}{|\Psi(\mu
  j)|^{q+2}} = \mu\sum_{j\in\ZZ}{f(j)} = \mu\int_{\RR}|\Psi(\mu
y)|^{q+2}dy + \\ &+& \mu^2\int_{\RR}\Psi(\mu y)|\Psi(\mu
y)|^{q}\Psi_x(\mu y)P_1(y)dy =\\ &=& \int_{\RR}|\Psi(x)|^{q+2}dx +\mu
\int_{\RR}\Psi(x)|\Psi(x)|^{q}\Psi_x(x)P_1(x/\mu)dx.
\end{eqnarray*}
Hence
\begin{displaymath}
R_H = \mu \int_{\RR}\Psi|\Psi|^{q}\Psi_x P_1(x/\mu)dx.
\end{displaymath}
A direct computation of the firsy and second differential shows that
\begin{eqnarray*}
dR_H(\Psi)[h] &=& \mu \int_{\RR}(\Psi|\Psi|^{q})_x h P_1(x/\mu)dx +
\mu \int_{\RR}\Psi|\Psi|^{q} h_x P_1(x/\mu)dx,\\
d^2R_H(\Psi)[h,h] &=& \mu \int_{\RR}(|\Psi|^{q})_x h^2 P_1(x/\mu)dx +
\mu \int_{\RR}|\Psi|^{q} (h^2)_x P_1(x/\mu)dx.\\
\end{eqnarray*}
The smallness is represented by the prefactor $\mu$: so Sobolev
embedding Theorems and $|P_1(x/\mu)|\leq 1/2$ yield
\eqref{e.p.A}. \qed

\bigskip

We have thus verified the assumptions of Theorem \ref{t.a.1} which
implies the existence and the estimate of the ST-mode for the case
$n=1$. The same statement for the P-mode follows by a translation of
the basis of $\E_\mu$
\begin{equation}
\label{Pmode-1}
\Psi(x) := \sum_{j\in\ZZ}{\psi_j}s_j(x/\mu + 1/2),\qquad\qquad
\psi_{-j+1} = \psi_j,
\end{equation}
with $s_j$ defined in \eqref{fe.1}.

\bigskip


\subsection{The case $n=2$}

\begin{figure}[t]
\begin{center}
\includegraphics[width=.6 \textwidth]{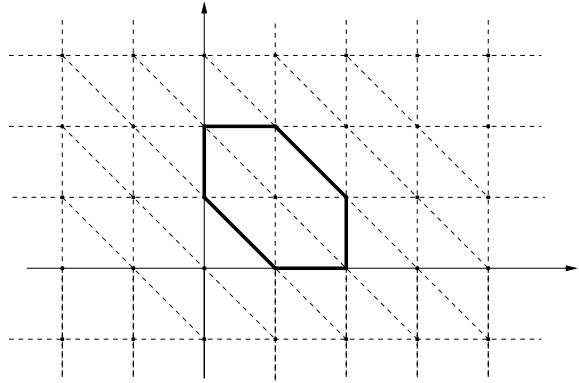}
\end{center}
\caption{Triangulation and finite element in the bidimensional lattice.}
\label{fig2}
\end{figure}

Let us take $\psi_{j,k}\in \Q_\mu$. For each multindex ${\bf l} =
(j,k)$, let us consider the function $s_{j,k}(x,y)$ which represents
the exagonal pyramid of height one centered in $(j,k)$ whose support
is the union of the six triangles of figure \ref{fig2}. More precisely
we define $T^+_{j,k}$ the triangle whose vertexes are
$(j,k),(j+1,k),(j,k+1)$ and $T^-_{j,k}$ the one whose vertexes are
$(j,k),(j-1,k),(j,k-1)$. Hence, for example, on $T^+_{j,k}$ the
function $s_{j,k}$ represents the plane in ${\RR}^3$
\begin{displaymath}
s_{j,k}(x,y) = -{x} -{y}+j+k+1.
\end{displaymath}
The set of functions $\{s_{j,k}(x/\mu,y/\mu)\}_{(j,k)\in\ZZ^2}$ is a
basis which generates a piecewise linear function $\Psi(x,y)$
interpolating $\psi_{j,k}$
\begin{equation}
\label{e.int.2}
\Psi(x,y) := \sum_{(j,k)\in\ZZ^2}{\psi_{j,k}s_{j,k}(x/\mu,y/\mu)}.
\end{equation}
Notice that on the triangle $T^{\pm}_{j,k}$ the function $\Psi$ is the
plane
\begin{equation}
\label{e.T+-}
\Psi(x,y) = \psi_{j,k} + (x-\mu j)\Psi_x +(y-\mu k)\Psi_y.
\end{equation}

\begin{definition}
\label{linspace.2}
We denote by $\E_\mu$ the linear space composed by the functions of
the form \eqref{e.int.2} with $\psi_{j,k}\in \Q_\mu$.
\end{definition}

The following Lemma gives the equivalence between {\bf the function
  space} $\E_\mu\subset H^1$ and {\bf the sequence space} $\Q_\mu$.

\begin{lemma}
\label{l.grad.2}
Let $\Psi\in\E_\mu$ then it holds true
\begin{equation}
\label{e.norm.der.2}
\int_{\RR^2}\Psi_x^2 = \sum_{(j,k)\in\ZZ^2}{(\psi_{j+1,k} -
  \psi_{j,k})^2},\qquad \int_{\RR^2}\Psi_y^2 =
\sum_{(j,k)\in\ZZ^2}{(\psi_{j,k+1} - \psi_{j,k})^2}.
\end{equation}
Moreover
\begin{equation}
\label{e.norm.ell2}
\mu^2\sum_{j,k}{\psi_{j,k}^2} = \int_{\RR^2} \Psi^2 +
\frac{\mu^2}6\int_{\RR^2}\tond{\Psi_x^2+\Psi_y^2 - \Psi_x \Psi_y}.
\end{equation}
\end{lemma}

\proof from \eqref{e.T+-} we have that on each triangle
$T^{\pm}_{j,k}$ it holds
\begin{displaymath}
\Psi_x = \pm\frac{(\psi_{j\pm1,k} - \psi_{j,k})}\mu,\qquad
\Psi_y = \pm\frac{(\psi_{j,k\pm1} - \psi_{j,k})}\mu.
\end{displaymath}
Formula \eqref{e.norm.ell2} follows from a direct computation as in
Lemma \ref{l.grad.1}.
\qed

\bigskip
\bigskip
The next three Lemmas provide the proof of the following main
\begin{proposition}
\label{p.2d}
Let $\Psi\in\E_\mu$ be as in \eqref{e.int.2} and let us define
\begin{eqnarray*}
G_c(\Psi) &:=& \int_{\RR^2}|\Psi|^{q+2},\\
G_d(\Psi) &:=& \mu^2\sum_{j,k\in\ZZ^2}{|\psi_{j,k}|^{q+2}},\\
R_G(\Psi) &:=&  G_c(\Psi) - G_d(\Psi);
\end{eqnarray*}
if $q\geq 1$ then $R_G\in\C^2(\E_\mu)$ and in any open set
$\U\subset\E_\mu$ one has
\begin{displaymath}
\norma{R_G}_{\C^2(\U)}\leq C_{\U} \mu.
\end{displaymath}
\end{proposition}

\bigskip
\bigskip

%
%

\begin{lemma}
\label{l.pot.2d}
Under the assumptions of Proposition \ref{p.2d} one has
\begin{displaymath} 
|R_G(\Psi)|\leq C \mu\norma{\Psi}_{H^1}^{q+2}.
\end{displaymath}
\end{lemma}

\proof Let us set $f(x,y)=|\Psi(x,y)|^{q+2}$ and let us take $(x,y)\in
T^\pm_{j,k}$, then we can use a Taylor expansion with integral
remainder
\begin{equation}
\label{a.1}
f(x,y) = f(\mu j,\mu k) + (x-\mu j)\int_0^1 f_x(\gamma(t,x,y))dt +
(y-\mu k)\int_0^1 f_y(\gamma(t,x,y))dt,
\end{equation}
where 
\[
\gamma(t,x,y) = (tx+(1-t)\mu j, ty+(1-t)\mu k) \qquad\qquad t\in[0,1]
\]
is the segment connecting $(\mu j,\mu k)$ with $(x,y)$ and lies in the
triangle $T^\pm_{j,k}$. Hence
\begin{eqnarray}
\int_{T^\pm_{j,k}}f(x,y) &=& \frac{\mu^2}2 f(\mu j,\mu k) +\nonumber\\
\label{parx}
&+&\int_{T^\pm_{j,k}}(x-\mu j)\int_0^1 f_x(\gamma(t,x,y))dtdxdy +\\
\label{pary} &+&\int_{T^\pm_{j,k}}(y-\mu k)\int_0^1 f_y(\gamma(t,x,y))dtdxdy.
\end{eqnarray}
By the initial definition of $f(x,y)$ one has
\begin{eqnarray}
\label{der_f}
|\partial_x f| &=& (q+2)|\Psi|^{q+1} |\Psi_x| = (q+2)|\Psi|^{q+1}
\frac{|\psi_{j\pm1,k}-\psi_{j,k}|}\mu,\\\nonumber |\partial_y f| &=&
(q+2)|\Psi|^{q+1} |\Psi_y| = (q+2)|\Psi|^{q+1}
\frac{|\psi_{j,k\pm1}-\psi_{j,k}|}\mu.\\
\end{eqnarray}
Since
\[
\int_{\RR^2}f(x,y) - \mu^2\sum_{j,k\in\ZZ^2}f(j,k) =
\sum_{j,k\in\ZZ^2}\tond{\int_{T^\pm_{j,k}}f(x,y)-\frac{\mu^2}2 f(j,k)}
\]
we can use \eqref{parx} and \eqref{pary} to estimate
\begin{eqnarray*}
|\int_{\RR^2}f(x,y) &-& \mu^2\sum_{j,k\in\ZZ^2}f(j,k)| \leq
\sum_{j,k\in\ZZ^2}|\int_{T^\pm_{j,k}}f(x,y)-\frac{\mu^2}2
f(j,k)|\leq\\ &\leq& \mu\sum_{j,k\in\ZZ^2}\int_{T^\pm_{j,k}}\int_0^1
|f_x(\gamma(t,x,y))| + |f_y(\gamma(t,x,y))|dtdxdy.\\ 
\end{eqnarray*}
By inserting \eqref{der_f} and observing that
\begin{equation}
\label{e.Psi.gamma}
\Psi\circ\gamma = \Psi(\gamma(t,x,y)) = \psi_{j,k} + t(x-\mu j)\Psi_x
+ t(y-\mu k)\Psi_y
\end{equation}
one may remove the integration along the segment
\[
\int_{T^\pm_{j,k}}\int_0^1 |f_x(\gamma(t,x,y))|\leq
\frac{|\psi_{j+1,k} -
  \psi_{j,k}|}\mu\int_{T^\pm_{j,k}}\quadr{|\psi_{j,k}| +
  |\langle\nabla\Psi,\delta_{j,k}\rangle|}^{q+1}
\]
where we set $\delta_{j,k}(x,y):= (x-\mu j,y-\mu k)$ and 
\[
\langle\nabla\Psi,\delta_{j,k}\rangle = (x-\mu j)\Psi_x + (y-\mu k)\Psi_y.
\]
Notice that from \eqref{e.T+-} one has
\[
|\psi_{j,k}|\leq |\Psi(x,y)| + \mu\tond{|\Psi_x|+|\Psi_y|}
\]
thus it is possible to estimate the argument of the integral
as follows \footnote{We here remind the inequalities for $a,b>0$
\[
(a+b)^s\leq 2^s(a^s+b^s),\qquad 0<s<1,
\]
and
\[
(a+b)^s\leq 2^{s-1}(a^s+b^s),\qquad s\geq 1;
\]
the second follows easily from the convexity of the function $g(x) =
x^s,\,x\in\RR^+$. The first is a direct consequence:
\[
(a+b)^s = \frac{(a+b)^{s+1}}{(a+b)}\leq
2^{s}\tond{\frac{a^{s+1}}{(a+b)}+\frac{b^{s+1}}{(a+b)}}\leq
2^{s}(a^s+b^s).
\]
}
\[
\quadr{|\psi_{j,k}| +
  |\langle\nabla\Psi,\delta_{j,k}\rangle|}^{q+1}\leq
4^q\quadr{|\Psi(x,y)|^{q+1} +
  \mu^{q+1}\tond{|\Psi_x|+|\Psi_y|}^{q+1}};
\]
hence recalling also that $|\Psi_x|+|\Psi_y|\leq\sqrt{2}|\nabla\Psi|$
one has
\begin{eqnarray*}
&&\mu\sum_{j,k\in\ZZ^2}\int_{T^\pm_{j,k}}\int_0^1 |f_x(\gamma(t,x,y))|
  + |f_y(\gamma(t,x,y))|dtdxdy=\\ &=&\mu
  C(q)\sum_{j,k\in\ZZ^2}\int_{T^\pm_{j,k}}|\nabla
  \Psi|\quadr{|\Psi(x,y)|^{q+1} + \mu^{q+1}|\nabla\Psi|^{q+1}}\leq
  \\ &\leq& \mu C(q)\sum_{j,k\in\ZZ^2}\quadr{\int_{T^\pm_{j,k}}|\nabla
    \Psi||\Psi(x,y)|^{q+1} + \mu^{q+1}\int_{T^\pm_{j,k}}|\nabla
    \Psi|^{q+2}}.
\end{eqnarray*}
The first sum gives
\begin{eqnarray*}
\sum_{j,k\in\ZZ^2}\int_{T^\pm_{j,k}}|\nabla \Psi||\Psi(x,y)|^{q+1} &=&
\int_{\RR^2}|\nabla \Psi||\Psi(x,y)|^{q+1}\leq\\&\leq&
\norma{\nabla\Psi}_{L^2}\norma{\Psi}_{L^{2q+2}(\RR^2)}^{q+1}
\leq\\ &\leq&C_1\norma{\nabla\Psi}_{L^2}\norma{\Psi}_{H^1}^{q+1}.
\end{eqnarray*}
Using Lemma \ref{l.grad.2} the second instead gives
\begin{eqnarray*}
\mu^{q+1}\sum_{j,k\in\ZZ^2}\int_{T^\pm_{j,k}}|\nabla \Psi|^{q+2}&=&
\mu^{q+3}\sum_{j,k\in\ZZ^2}|\nabla \Psi|^{q+2}
=\\ &=&\mu\sum_{j,k\in\ZZ^2}\quadr{(\psi_{j+1,k} -
  \psi_{j,k})^2+(\psi_{j,k+1} -
  \psi_{j,k})^2}^{1+\frac{q}2}\leq\\ &\leq&
C_2\mu\norma{\nabla\Psi}_{L^2}^{q+2}.
\end{eqnarray*}
Collecting the above estimates we obtain
\[
\big|\int_{\RR^2}|\Psi(x,y)|^{q+2} -
\mu^2\sum_{j,k\in\ZZ^2}|\psi_{j,k}|^{q+2}\big| \leq \mu C_1
\norma{\nabla\Psi}_{L^2}\norma{\Psi}_{H^1}^{q+1} + \mu^2 C_2
\norma{\nabla\Psi}_{L^2}^{q+2}
\]
which finally gives
\[
\big|\int_{\RR^2}|\Psi(x,y)|^{q+2} -
\mu^2\sum_{j,k\in\ZZ^2}|\psi_{j,k}|^{q+2}\big|\leq \mu
C\norma{\Psi}_{H^1}^{q+2}.
\]
\qed

\bigskip
\bigskip

\begin{lemma}
\label{l.pot.2d.C1}
Under the assumptions of Proposition \ref{p.2d} one has also
\begin{equation}
\label{e.pot.2d.C1}
\norma{R_G'(\Psi)}_{\L(\E_\mu,\RR)}\leq C\mu\norma{\Psi}_{H^1}^{q+1}.
\end{equation}
\end{lemma}

\proof A direct computation easily gives for any $h\in\E_\mu$
\[
G_c'(\Psi)[h] = \int_{\RR^2}|\Psi|^q\Psi h,\qquad\qquad G_d'(\Psi)[h] =
\mu^2\sum_{j,k}|\psi_{j,k}|^q\psi_{j,k} h_{j,k}
\]
with obviously $h_{j,k} = h(\mu j,\mu k)$. In order to estimate
\[
\norma{G_c'(\Psi) - G_d'(\Psi)}_{\L(\E_\mu,\RR)} =
\sup_{h\not=0}\frac{|G_c'(\Psi)[h] - G_d'(\Psi)[h]|}{\norma{h}_{H^1}}
\]
we need to control
\begin{equation}
\label{e.dRG}
|G_c'(\Psi)[h] - G_d'(\Psi)[h]| = |\int_{\RR^2}|\Psi|^q\Psi h -
\mu^2\sum_{j,k}|\psi_{j,k}|^q\psi_{j,k} h_{j,k}|.
\end{equation}
We proceed in the same way as in Lemma \ref{l.pot.2d}, exploiting the
fact that also $h(x,y)\in\E_\mu$. We thus define in this case $f(x,y)
= |\Psi(x,y)|^q\Psi(x,y)h(x,y)$, so that
\[
|f_x|+|f_y|\leq q|\Psi|^q(|\Psi_x|+|\Psi_y|)|h| + |\Psi|^{q+1}|\nabla
h|.
\]
We recall that
\begin{eqnarray*}
|\Psi\circ\gamma| &\leq& |\Psi|+2\mu(|\Psi_x|+|\Psi_y|),\\
|h\circ\gamma| &\leq& |h|+2\mu|\nabla h|,
\end{eqnarray*}
hence \eqref{e.dRG} can be split into four terms
\begin{eqnarray*}
\sum_{\ZZ^n}\int_{T^{\pm}_{j,k}}\int_0^1{|f_x\circ\gamma|+|f_y\circ\gamma|}
&\leq& \int_{\RR^2}|\Psi|^q|\nabla\Psi|(|h|+\mu|\nabla h|) +\\ &+&
\mu^q\int_{\RR^2}|\nabla\Psi|^{q+1}(|h|+\mu|\nabla h|) +\\ &+&
\int_{\RR^2}|\Psi|^{q+1}|\nabla h|+\\ &+&
\mu^{q+1}\int_{\RR^2}|\nabla\Psi|^{q+1}|\nabla h|.
\end{eqnarray*}
\begin{enumerate}
\item using Schwarz and $(\nabla|\Psi|^{q+1})^2 =
  C(q)|\nabla\Psi|^2|\Psi|^{2q}$ we get
\begin{eqnarray*}
\int_{\RR^2}|\Psi|^q|\nabla\Psi|(|h|+\mu|\nabla h|) &\leq& C_1
\sqrt{\int_{\RR^2}|\Psi|^{2q}|\nabla\Psi|^2}\norma{h}_{H^1}\leq\\ &\leq&
C_2 \norma{\Psi}_{H^1}^{q+1}\norma{h}_{H^1};
\end{eqnarray*}

\item using Schwarz and $\ell^2\hookrightarrow\ell^s,\,s>2$
\begin{eqnarray*}
\mu^q\int_{\RR^2}|\nabla\Psi|^{q+1}(|h|+\mu|\nabla h|) &\leq&
C_1\mu^q\sqrt{\int_{\RR^2}|\nabla\Psi|^{2q+2}}\norma{h}_{H^1}\leq\\ &\leq&
C_2 \norma{\nabla\Psi}_{L^2}^{q+1}\norma{h}_{H^1};
\end{eqnarray*}
\item as in 1.
\begin{displaymath}
\int_{\RR^2}|\Psi|^{q+1}|\nabla h|\leq
\norma{\Psi}_{L^{2q+2}}^{q+1}\norma{\nabla h}_{L^2}\leq
C_1\norma{\Psi}_{H^1}^{q+1}\norma{h}_{H^1};
\end{displaymath}

\item as in 2.
\begin{eqnarray*}
\mu^{q+1}\int_{\RR^2}|\nabla\Psi|^{q+1}|\nabla h| &\leq&
\mu^{q+1}\sqrt{\int_{\RR^2}|\nabla\Psi|^{2q+2}}\norma{\nabla
  h}_{L^2}\leq\\ &\leq&
C_1\norma{\nabla\Psi}_{L^2}^{q+1}\norma{h}_{H^1}.
\end{eqnarray*}
\end{enumerate}
Collecting we get
\[
|R_G'(\Psi)[h]|\leq C(q)\mu\norma{\Psi}_{H^1}^{q+1}\norma{h}_{H^1},
\]
hence the thesis.
\qed

\bigskip
\bigskip

\begin{lemma}
\label{l.pot.2d.C2}
Under the assumptions of Propositions \ref{p.2d} one has also
\begin{equation}
\label{e.pot.2d.C2}
\norma{R_G''(\Psi)}_{\L^2(\E_\mu,\RR)}\leq C\mu\norma{\Psi}_{H^1}^{q}.
\end{equation}
\end{lemma}

\proof Also in this case, a direct computation easily gives for any
$h\in\E_\mu$
\[
G_c''(\Psi)[h,h] = \int_{\RR^2}|\Psi|^q h^2,\qquad\qquad
G_d''(\Psi)[h,h] = \mu^2\sum_{j,k}|\psi_{j,k}|^q h_{j,k}^2.
\]
In order to estimate
\[
\norma{G_c''(\Psi) - G_d''(\Psi)}_{\L^2(\E_\mu,\RR)} =
\sup_{h\not=0}\frac{|G_c''(\Psi)[h,h] -
  G_d''(\Psi)[h,h]|}{\norma{h}^2_{H^1}}
\]
we need to control
\begin{equation}
\label{e.d2RG}
|G_c''(\Psi)[h,h] - G_d''(\Psi)[h,h]| = |\int_{\RR^2}|\Psi|^q h^2 -
\mu^2\sum_{j,k}|\psi_{j,k}|^q h_{j,k}^2|.
\end{equation}
We proceed as in the previous Lemmas, by defining 
\begin{displaymath}
f(x,y) := |\Psi(x,y)|^qh^2(x,y),
\end{displaymath}
so that
\[
|f_x|+|f_y|\leq q|\Psi|^{q-1}|\nabla\Psi|h^2 +
2|\Psi|^{q}|h\nabla h|.
\]
We distinguish the case $q=3$, which is easier, and $q>3$.
\begin{description}
\item[$q=1$]
In this case we have
\begin{eqnarray*}
\sum_{\ZZ^n}\int_{T^{\pm}_{j,k}}\int_0^1{|f_x\circ\gamma|+|f_y\circ\gamma|}
&\leq& 2\int_{\RR^2}|\nabla\Psi||h|^2 +\\ &+&
4\mu^2\int_{\RR^2}|\nabla\Psi||\nabla h|^2 +\\ &+&
\int_{\RR^2}|\Psi||\nabla h|(|h|+\mu|\nabla h|)+\\ &+&
\mu\int_{\RR^2}|\nabla\Psi||\nabla h|(|h|+\mu|\nabla h|).
\end{eqnarray*}
The thesis can be obtained using Schwarz and observing that
\[
\int_{\RR^2}|\nabla h|^4 \leq
\frac{c_1}{\mu^2}\sum_{j,k}\quadr{(h_{j+1,k}-h_{j,k})^2+(h_{j,k+1}-h_{j,k})^2}^2\leq
\frac{c_2}{\mu^2}\norma{\nabla h}_{L^2}^4.
\]
\item[$q>1$] The steps are the same as usual; the only difference is
  that we have to deal with
\[
\int_{\RR^2}|\Psi|^{2\sigma}|\nabla\Psi|^2,\qquad\qquad \sigma>0,
\]
but it is enough to notice again that the integral above is the
(square) $L^2$ norm of $\nabla|\Psi|^{1+\sigma}$, thus
\[
\int_{\RR^2}|\Psi|^{2\sigma}|\nabla\Psi|^2\leq
C\norma{\Psi}_{H^1}^{2+2\sigma}.
\]
\end{description}
\qed

\bigskip
\bigskip

This concludes the case related to the construction and approximation
of the ST-mode. The other three modes (the P-mode and the two H-mode)
are obtained by translation of the basis $s_{j,k}$ either in one or in
both the two directions.

\bigskip
\bigskip


\subsection{Proof of Theorem \ref{main.theo}.}
We begin with the following

\begin{definition}
\label{d.ere}
Let $n=1,2$ and consider $\psi\in H^2(\RR^n)\hookrightarrow\C^0$ on
$\E_\mu$. We define
\begin{equation}
\label{e.proj.1}
\Pi_\mu:\psi \mapsto \Pi_\mu \psi = \sum_{{\bf l}\in\ZZ^n}
\psi(\mu {\bf l}) s_{\bf l}(x/\mu),\qquad x\in\RR^n
\end{equation}
the projection of $H^2(\RR^n)$ on $\E_\mu$. By classical results on
polynomial approximation in Sobolev spaces (Chapter~4 of
\cite{BrenScot}) one has
\begin{equation}
\label{e.approx}
\norma{\Pi_\mu \psi - \psi}_{H^1} \leq C\mu\norma{\psi}_{H^2}.
\end{equation}
\end{definition}

We need also a simple lemma to obtain the second estimate of
\eqref{est}

\begin{lemma}
\label{l.sup}
For any ${\bf l}\in\ZZ^n$ we have
\begin{equation}
\label{e.sup}
\left|\psi_{\bf l}\right|\leq 2\mu^{\frac12-\frac{n}2} \nqmu{\psi}.
\end{equation}
\end{lemma} 

\proof We write the proof for the case $n=2$. The case $n=1$ is
simpler. Denote ${\bf l}=(j,k)$; one has
\begin{displaymath}
\psi_{j,k}^2 = \sum_{m=-\infty}^h(\psi_{m,k}^2 - \psi_{m-1,k}^2) =
\sum_{m=-\infty}^h(\psi_{m,k} - \psi_{m-1,k})(\psi_{m,k} + \psi_{m-1,k})
\end{displaymath}
which gives
\begin{eqnarray*}
\sup_{(j,k)\in\ZZ^2}\psi_{j,k}^2 &\leq& 4 \sqrt{\sum_{m\in\ZZ}\psi_{m,k}^2}
\sqrt{\sum_{m\in\ZZ}(\psi_{m+1,k} - \psi_{m,k})^2} \leq
\\
&\leq& \tond{2\sqrt{\norma
  {\psi}_{\ell^2}}\mu^{\frac{1}{2}-\frac{n}{4}} \left[ \mu^n
  \frac{\left\langle \psi;-\Delta\psi  \right\rangle_{\ell^2 }}{\mu^2}
  \right]^{1/4}}^2 
\end{eqnarray*}
and \eqref{e.sup}.\qed

\bigskip
\bigskip

Now we easily verify the hipothesis of the abstract Theorem
\ref{t.a.1}. First, we define $\psi_c$ as the (smooth) solution of
\eqref{lagrange.c} and $\psi_0 = \nu\Pi_\mu\psi_c$, with $\nu$ such
that $N_{0,\mu}(\psi_0)=1$. Then, condition \eqref{coerc} follows from
Lemma \ref{NDegen} while condition \eqref{e.a.2} comes from the above
\eqref{e.approx}. Finally, requirement {ii} is given by Lemmas
\ref{l.grad.1} and \ref{l.grad.2} and by Propositions \ref{p.1d} and
\ref{p.2d}. This directly gives the first of \eqref{est}. The second
of \eqref{est} is a byproduct of either the first and Lemma
\ref{l.sup}, indeed
\begin{displaymath}
\left|{\psi^i_{\bf l}(t) - \Psi^i_{\bf l}(t)} \right|\leq
2\mu^{\frac12-\frac{n}2}\nqmu{\psi^i-\Psi^i}\leq C
\mu^{\frac32-\frac{n}2}.
\end{displaymath}
\qed

\bigskip\bigskip

\hrule

\bigskip\bigskip

{\bf Acknowledgments:} we warmly thank Simone Paleari for his constant
and stimulating comments, and Livio Pizzocchero and Massimo Tarallo
for their valuable help. We also thank P. Kevrekidis for his
bibliographic informations. This work has been partially supported by
PRIN 2007B3RBEY ``Dynamical Systems and applications''.

\bigskip\bigskip\bigskip

\hrule

\bigskip\bigskip\bigskip



\begin{thebibliography}{CJK{\etalchar{+}}08}

\bibitem[AKK01]{AubKopKad01}
S.~Aubry, G.~Kopidakis, and V.~Kadelburg, \emph{Variational proof for hard
  discrete breathers in some classes of {H}amiltonian dynamical systems},
  Discrete Contin. Dyn. Syst. Ser. B \textbf{1} (2001), no.~3, 271--298.

\bibitem[AP95]{AmbPro}
A.~Ambrosetti and G.~Prodi, \emph{A primer of nonlinear analysis}, Cambridge
  Studies in Advanced Mathematics, vol.~34, Cambridge University Press,
  Cambridge, 1995, Corrected reprint of the 1993 original.

\bibitem[BCP02]{BCP}
D.~Bambusi, A.~Carati, and A.~Ponno, \emph{The nonlinear {S}chr\"odinger
  equation as a resonant normal form}, DCDS-B \textbf{2} (2002), 109--128.

\bibitem[BCP09]{BCP09}
D.~Bambusi, A.~Carati, and T.~Penati, \emph{Boundary effects on the dynamics of
  chains of coupled oscillators}, Nonlinearity \textbf{22} (2009), 923--946.

\bibitem[BL83]{BeLio83}
H.~Berestycki and P.-L. Lions, \emph{Nonlinear scalar field equations. {I}.
  {E}xistence of a ground state}, Arch. Rational Mech. Anal. \textbf{82}
  (1983), no.~4, 313--345.

\bibitem[BLP81]{BeLioPel81}
H.~Berestycki, P.-L. Lions, and L.~A. Peletier, \emph{An {ODE} approach to the
  existence of positive solutions for semilinear problems in {${\bf
  R}\sp{N}$}}, Indiana Univ. Math. J. \textbf{30} (1981), no.~1, 141--157.

\bibitem[BP06]{BP06}
D.~Bambusi and A.~Ponno, \emph{On metastability in {FPU}}, Comm. Math. Phys.
  \textbf{264} (2006), no.~2, 539--561.

\bibitem[BS08]{BrenScot}
S.~C. Brenner and L.~R. Scott, \emph{The mathematical theory of finite element
  methods}, third ed., Texts in Applied Mathematics, vol.~15, Springer, New
  York, 2008.

\bibitem[CGM78]{ColGlaMar78}
S.~Coleman, V.~Glaser, and A.~Martin, \emph{Action minima among solutions to a
  class of {E}uclidean scalar field equations}, Comm. Math. Phys. \textbf{58}
  (1978), no.~2, 211--221.

\bibitem[CJK{\etalchar{+}}08]{JmsKevr08}
J.~Cuevas, G.~James, P.~G. Kevrekidis, B.~A. Malomed, and B.~Sanchez-Rey,
  \emph{Approximation of solitons in the discrete {NLS} equation}, J. Nonlinear
  Math. Phys. \textbf{15} (2008), no.~suppl. 3, 124--136.

\bibitem[FGJS04]{FGJS04}
J.~Fr{\"o}hlich, S.~Gustafson, B.~L.~G. Jonsson, and I.~M. Sigal,
  \emph{Solitary wave dynamics in an external potential}, Comm. Math. Phys.
  \textbf{250} (2004), no.~3, 613--642.

\bibitem[FP99]{FPI}
G.~Friesecke and R.~L. Pego, \emph{Solitary waves on {F}ermi-{P}asta-{U}lam
  lattices. {I}. {Q}ualitative properties, renormalization and continuum
  limit.}, Nonlinearity \textbf{12} (1999), 1601--1627.

\bibitem[FP02]{FPII}
\bysame, \emph{Solitary waves on {F}ermi-{P}asta-{U}lam lattices. {II}.
  {L}inear implies nonlinear stability.}, Nonlinearity \textbf{15} (2002),
  1343--1359.

\bibitem[FP04a]{FPIII}
\bysame, \emph{Solitary waves on {F}ermi-{P}asta-{U}lam lattices. {III}.
  {H}owland-type {F}loquet theory.}, Nonlinearity \textbf{17} (2004), 207--227.

\bibitem[FP04b]{FPIV}
\bysame, \emph{Solitary waves on {F}ermi-{P}asta-{U}lam lattices. {IV}. {P}roof
  of stability at low energy.}, Nonlinearity \textbf{17} (2004), 229--251.

\bibitem[FW98]{FlaWil98}
S.~Flach and C.~R. Willis, \emph{Discrete breathers}, Phys. Rep. \textbf{295}
  (1998), no.~5, 181--264.

\bibitem[HW08]{HofWay08}
A.~Hoffman and C.~E. Wayne, \emph{Counter-propagating two-soliton solutions in
  the {F}ermi-{P}asta-{U}lam lattice}, Nonlinearity \textbf{21} (2008), no.~12,
  2911--2947.

\bibitem[Jam03]{Jam03}
G.~James, \emph{Centre manifold reduction for quasilinear discrete systems}, J.
  Nonlinear Sci. \textbf{13} (2003), no.~1, 27--63.

\bibitem[KRB01]{KevrRasBis01B}
P.~G. Kevrekidis, K.~O. Rasmussen, and A.~R. Bishop, \emph{The discrete
  nonlinear schr\"odinger equation: a survey of recent results}, Inernational
  Journal of Modern Physics B \textbf{15} (2001), no.~21, 2833--2900.

\bibitem[KSM92]{KSM92}
P.~Kirrmann, G.~Schneider, and A.~Mielke, \emph{The validity of modulation
  equations for extended systems with cubic nonlinearities}, Proc. Roy. Soc.
  Edinburgh Sect. A \textbf{122} (1992), no.~1-2, 85--91.

\bibitem[MP08]{MizPeg08}
T.~Mizumachi and R.~L. Pego, \emph{Asymptotic stability of {T}oda lattice
  solitons}, Nonlinearity \textbf{21} (2008), no.~9, 2099--2111.

\bibitem[Sch98]{Sch98}
G.~Schneider, \emph{Justification of modulation equations for hyperbolic
  systems via normal forms}, NoDEA Nonlinear Differential Equations Appl.
  \textbf{5} (1998), no.~1, 69--82.

\bibitem[SW00]{SW99}
G.~Schneider and C.~E. Wayne, \emph{Counter-propagating waves on fluid surfaces
  and the continuun limit of the {F}ermi {P}asta {U}lam model}, Proceedings of
  the International Conference on Differential Equations, Berlin 1999 (Rivere
  Edge, NJ), World Scientific, 2000.

\bibitem[Wei99]{Wei_MI99}
M.~I. Weinstein, \emph{Excitation thresholds for nonlinear localized modes on
  lattices}, Nonlinearity \textbf{12} (1999), no.~3, 673--691.

\end{thebibliography}

\newcommand{\etalchar}[1]{$^{#1}$}
\providecommand{\bysame}{\leavevmode\hbox to3em{\hrulefill}\thinspace}
\providecommand{\MR}{\relax\ifhmode\unskip\space\fi MR }
\providecommand{\MRhref}[2]{%
  \href{http://www.ams.org/mathscinet-getitem?mr=#1}{#2}
}
\providecommand{\href}[2]{#2}

\end{document}